\documentclass[review,3p]{elsarticle}

\usepackage{amsmath}
\usepackage{amssymb,amsfonts,textcomp}
\usepackage{paralist}

\numberwithin{equation}{section}

\begin{document}

\begin{frontmatter}

\title{Stabilization with residual-free bubbles for advection-dominated transport equations}

\author[uva]{I. Kryven\corref{cor1}}
\ead{i.kryven@uva.nl}
\author[lnu]{V. Kukharskyy}
\ead{vetaley@franko.lviv.ua}
\author[lnu]{Ya. Savula}
\ead{savula@lnu.edu.ua} 
\cortext[cor1]{Corresponding author} 

\address[lnu]{Lviv National University, Universytetska 1, 79000 Lviv, Ukraine}
\address[uva]{University of Amsterdam, Science Park 904, 1098 XH, Amsterdam, The Netherlands}

\begin{abstract}
An improved numerical scheme is proposed for advection-dominated advection-diffusion problem. The scheme is based on Galerkin  finite element method (FEM) with basis enriched with approximations to residual-free bubbles. The stabilisation effect of the numerical scheme was studied on several benchmark problems with high-P\'eclet numbers.  Comparison to traditional hp-FEM, reveals improved stability and computational performance.
\end{abstract}

\end{frontmatter}
%=============================Section One==========================================
\section{Introduction}
This paper proposes a numerical approximation method for solving a steady linear advection-diffusion  problem, described by the partial differential equation (PDE),
\begin{equation}\label{Lu}
\begin{array}{l}
Lu:=-k\Delta u+\mathbf{w}\cdot \nabla u=f$ in $\Omega \\
u|_{u \in \partial\Omega}=0.
\end{array}
\end{equation}
Here the unknown $u$ is a real function defined on a bounded domain $\Omega \in \mathbb{R}^{2}$ with Lipschitz continuous boundary $\partial \Omega $; $f \in L^{2}(\Omega)$ is the source term; the diffusion coefficient $k>0$ and the vector-valued advection filed $\mathbf{w}$ are assumed to be piecewise constant functions on $\Omega$.
The equation \eqref{Lu} becomes a singularly perturbed problem when $k \ll |\mathbf{w}|h$, where $h$ is a discretization step of a numerical method. It's known that in this case the standard Galerkin FEM based on linear or bilinear polynomials doesn't solve the problem adequately: non-physical oscillations are exhibited \cite{KuharskySavulaHolovachEN}.

A family of enriched methods has been developed in order to overcome such artefacts. We will address the most efficient methods of this class: the hp-FEM and the residual-free bubble (RFB) method. 
The hp-FEM is a generalization of FEM based on piecewise-polynomial approximations, which employs elements of variable size $h$ and requires a special system of basis functions of order $p$ \cite{BabuskaGuo}. 

	In the RFB method the basis of the standard FEM is enriched with so called 'bubble' functions, that satisfy exactly or approximately an equation involving the differential operator \eqref{Lu} in the interior of each element and vanish on their boundaries \cite{BrezziRusso}. Such approach requires solving a PDE sub-problem for every element of FEM to define the enrichment. The RFB method is often used for advection dominated  advection-diffusion problems \cite{BrezziMariniRusso,FrancaHawang,AsensioRussoSangalli,AyusoMarini}. 

The accuracy of approximation for RBF sub-problems strongly correlates with convergence speed of the enriched method.
Indeed, if for every element $K$ we set the bubble function $u^b_K \equiv 0$, the RFB method will be reduced to the standard FEM. On the other hand, if we exactly satisfy the sub-problem's differential equation, the 
RFB method gives an almost optimal approximation in $H^1$ \cite{BarboneHarari}. For this reason, the strategy of sub-problem approximation for the RFB method is very important.

In some cases the standard FEM (often with a refined mesh) is employed for solving RBF sub-problems on every element \cite{FrancaRamalhoValentin,Sangalli,DolbowFranca}. Such approach exhibits stable solutions for mesh P\'eclet numbers up to $10^5$.

Involving a bilinear approximation on sub-meshes for constructing the bubble enrichment $u^b_K$, brings certain difficulties:\\
\begin{inparaenum}[(i)]
\item using $u^b_K$ in the variational equation requires computing its gradients, but being a
piece-linear function, differentiation of $u^b_K$ brings well-known problems with stability; \\
\item for the same reason numerical integration of expressions with $u^b_K$ and its derivatives has to be done with a low order quadrature (i.e trapezoidal rule instead of Gauss-Legendre quadrature);\\  
\item finally, if the operator involved in \eqref{Lu} is strongly non-symmetric which is the main case addressed in this paper, the FEM needs to be applied with a small discretization step to avoid artefacts. This reflects dramatically on the overall  performance time, since the sub-problem needs to be solved on every element.\\
\end{inparaenum}
In this paper we will search for a more efficient approaches  to approximate RFB sub-problems.

\section{The RFB framework}
Let's denote a bilinear form $a:H^1_0\times H^1_0\to \mathbb{R}$
\begin{equation}\label{auv}
a(u,v):=<k\nabla u,\nabla v>+<\mathbf{w}\cdot \nabla u,v>,
\end{equation}
then the variational formulation for \eqref{Lu} is: to find $u^h\in H^1_0$, such that 
\begin{equation}\label{auv1}
a(u,v)=<f,v>, \ \forall v \in H^1_0(\Omega).
\end{equation}
We assume that in sense of measures $\nabla \cdot \mathbf{w}\le 0$ which guarantees the well posedness of \eqref{Lu} for any $k>0$: indeed it makes \eqref{auv} coercitive on $H^1_0(\Omega)$ {\cite{AsensioRussoSangalli}
\begin{equation}\label{coercitive_a}
a(v,v)\ge k|\mathbf{w}|_{H_{0}^{1}}^{2}.
\end{equation}
Let $\zeta^h=\{K\}$ be a standard partition of $\Omega$ in quadrilateral elements $K$, where $h=\underset{K in \zeta^h}{\mathit{max}}\mathit{diam}K$ -- mesh diameter. We assume $\zeta^h$ to be admissible (i.e. non-overlapping elements, their union reproduces domain, etc.) and shape regular (satisfying a minimal angle condition). We also choose a space $V^h\subset H_0^1(\Omega )$, related to the choice of partition $\zeta^h$. Thus the Galerkin method for the problem \eqref{Lu} is: to find $u^h\in V^h$, such that
\begin{equation}\label{DiscGalorkin}
a(u^{h},v^{h})=<f,v^{h}>\ \forall v^{h}\in V^{h}.
\end{equation}
We also assume that the diffusion coefficient and the advection vector field are constant on every $K$ from $\zeta^h$.
Let $S^1(K)$ be the space spanned by bilinear functions defined on the partition element $K$, $S^{1}(K)=\mathit{span}\{\phi_1^1,\phi _2^1,\phi _3^1,\phi _4^1\}$. Therefore an arbitrary function from $S^{1}(K)$ could be presented as
\begin{equation}\label{u1K}
u_K^1=\sum_{i=1}^4c_i\phi _i^1, \forall K\in
\zeta ^{h}.
\end{equation}
Let's split the approximation space $V^h$ into a composition of the finite dimension space $S^1(K)$ and the enrichment space of bubble functions
\begin{equation}\label{Vh}
	V^{h}(\Omega ) = V_{1}^{h}(\Omega )\oplus B^{h}(\Omega ),
\end{equation}
where
\begin{equation}\label{Vh1}
	\begin{array}{l}
		V_{1}^{h}(\Omega )=\{v\in H_{0}^{1}(\Omega):v|_{K}\in S^{1}(K), \forall K\in \zeta ^{h}\}, 
	\end{array}
\end{equation}
\begin{equation}\label{Bh}
	B^{h}(\Omega )=\{v\in H_{0}^{1}(\Omega):v|_{\partial K}=0,\forall K\in \zeta ^{h}\}.
\end{equation}
An element of \eqref{Vh} could be written as
\begin{equation}\label{uh}
V^{h}(\Omega )\ni u^{h}=u^{1}+u^{b}, 
\end{equation}
where $u^{1}\in V_{1}^{h}(\Omega ),\ u^{b}\in B^{h}(\Omega )$.
Therefore we can rewrite \eqref{DiscGalorkin} as
\begin{equation}\label{RBFvariation}
\begin{array}{l}
a(u^{h},v^{1}+v^{b})=<f,v^{1}+v^{b}>,\\ \forall v^{1}\in V_{1}^{h},\ v^{b}\in
B^{h},
\end{array}
\end{equation}
which leads to
\begin{equation}\label{RBFvariation1}
a(u^{h},v^{1})=<f,v^{1}>,\ \forall v^{1}\in V_{1}^{h}(\Omega);
\end{equation}
\begin{equation}\label{RBFvariation2}
a(u^{h},v^{b})=<f,v^{b}>,\ \forall v^{b}\in B^{h}(\Omega ).
\end{equation}
Taking \eqref{Bh} into account  we transform the variational formulation \eqref{RBFvariation2} to a set of PDE problems:
for every $K\in \zeta ^{h}$ find $u_{K}^{b}\in H^1_0(K)$, such that
\begin{equation}\label{LocalProblem}
\begin{array}{l}
Lu_{K}^{b}=f-Lu_{K}^{1} \text{ on } K, u_{K}^{1}\in S^{1}(K);\\
u_{K}^{b}|_{\partial K}=0.
\end{array}
\end{equation}
Let's denote the operator $F_K(x):S^{1}(K)\rightarrow H^1_0(K)$, that maps an element $u_{K}^{1}\in S^{1}(K)$ into 
a solution of the sub-problem \eqref{LocalProblem}. Evidently, $F_K$ is a linear non-invertible operator.

Let's denote $u^b_K$ as a linear combination of elements from $B^{h}(\Omega)$,
using the same coefficients as in \eqref{u1K} 
\begin{equation}\label{uRBFK}
u_K^b=\sum_{i=1}^4c_i\phi _i^b,\ \phi _i^b\in
H^1_0(K), \forall K\in \zeta^h.
\end{equation}
Using \eqref{uh}, \eqref{LocalProblem}, \eqref{u1K}, \eqref{uRBFK} and properties of operator $F_K$ we come up with
\begin{equation}\label{uRBFh}
\begin{array}{l}
u^{h}|_{K}=\sum\limits_{i=1}^{4}c_{i,j}(I+F_K)(\phi_i^1),
\\ i=1..4, K \in \zeta^h,

\end{array}
\end{equation}
here $I$ denotes the identity operator.

Generally, to apply the RFB method one needs to perform the following steps:\\
\begin{inparaenum}[(i)]
\item define partition $\zeta^h$;\\
\item calculate the values of $\phi _i^b=F_K(\phi
_{i}^{1}),\ i=1..4,$ for all $K\in \zeta ^{h}$;\\
\item using standard Galerkin method calculate $c_{i,j}$ from \eqref{RBFvariation1}, taking  the basis enrichment \eqref{uRBFh} into account.\\
\end{inparaenum}
We refer to \eqref{u1K}, \eqref{uRBFK} as the exact RFB method, because generally speaking, $\mathit{dim}(B^{h}(\Omega ))=\infty $. Even if we set degrees of freedom for the bubble space to be a finite number, we still need to compute $F_K(x)$ using a discrete method. The FEM on sub-mesh is commonly used for this purpose \cite{SangalliSmallScales,CariusMadureira,BrezziMariniRusso, AsensioRussoSangalli,FrancaRamalhoValentin}. As we will see later, using a spectral method often shows better stability for problems with high P\'eclet number.

\section{Spectral methods for solving RFB sub-problems}
With no loss of generality we consider $K$ to be
\begin{equation}\label{K}
\{(\xi ,\eta ):-1\le \xi ,\eta \le 1\}.
\end{equation}
We will use a discrete variation equation over a finite space instead of \eqref{RBFvariation2}.
\begin{equation}\label{aJ}
\begin{array}{l}
a(u^p,v^p)=<f,v^p>-<\mathbf{w}\cdot \nabla u_{K}^{1},v^p>,\\ \forall v^p \in B^p_{0}(K),
\end{array}
\end{equation}
where space $B^p_0$ is spanned by the hierarchy basis\cite{ShephardDeyFlaherty} 
\begin{equation}\label{Bs0}
B_{0}^p(K)=\mathit{span}\{M_{i,j}^r,\;r=1..p,\;i+j=r\},
\end{equation}
\begin{equation}\label{M}
\begin{array}{l}
M_{i,j}^p=(1-\xi ^2)(1-\eta ^2)L_i(\xi)L_j(\eta), \;i+j=p-1,
\end{array}
\end{equation}
$L_p(\tau)$ is the Legendre polynomial of $p^{th}$ order.
Ansatz for approximated  solution is:
\begin{equation}\label{UUs}
u^{p}=\sum_{r=1}^p\sum_{i+j=r-1}c_{i,j}^{r}M_{i,j}^{r},
\end{equation}
Thus, by fixing the approximation order $p$ and substituting \eqref{UUs} into variational equation \eqref{aJ}, \eqref{Bs0}, we arrive with a full discrete scheme for computing $c_{i,j}$:
\begin{equation}\label{SLAR}
\begin{array}{l}
\displaystyle\sum_{r=1}^{p}\displaystyle\sum_{i+j=r-1}c_{i,j}^{r}a(M_{i,j}^{r},M_{i_{1,}j_{1}}^{r_{1}})=\\=<f,M_{i_{1,}j_{1}}^{r_{1}}>-<\mathbf{w}\cdot
\nabla
u^{1},M_{i_{1,}j_{1}}^{r_{1}}>,\\ r_{1}=1...p,i_{1}+j_{1}=r_{1}-1.
\end{array}
\end{equation}
which is a linear system with $(p+p^2)/2$ equations, defined by a square, invertible and full matrix.

Proposed strategy has several benefits:\\
\begin{inparaenum}[(i)]
\item The approximation of the solution $u^p$ belongs to $C^p$ what makes an accurate gradient calculation possible;\\
\item the numerical solutions with high smoothness are known to exhibit better stability when approximating a solution of a PDE with strongly non-symmetric operators involved \cite{IdelsohnOnateCalvoPin};\\
\item it is relatively easy to apply an optimal, in respect to number and position of nodes,  integration quadrature to $u^p$ and to its first partial derivatives, as they are polynomials of degree $p$ and $p-1$ respectively;\\
\item there is no need to recalculate the matrix of the system \eqref{SLAR} for every finite element, but its right part column. Consequently, the computation time is saved.
\end{inparaenum}

\section{Numerical tests}
\subsection{A two-dimensional problem with symmetric boundary conditions}

\begin{figure}
\center
\includegraphics[width=70mm]{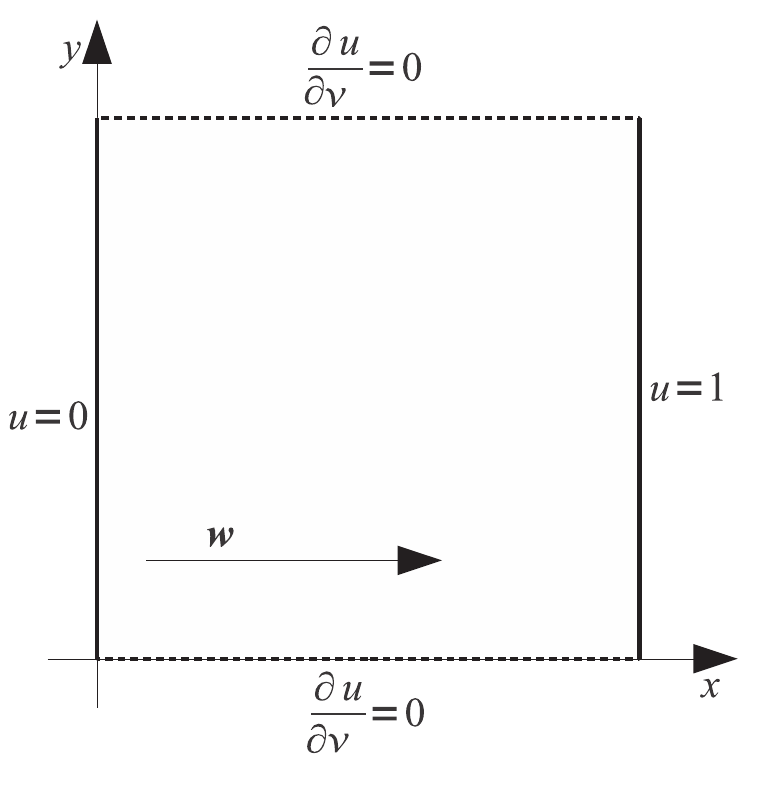}
\caption{2D domain $\Omega$ with symmetric boundary condition}
\label{omega0}
\end{figure}
For benchmark purposes, we provide simulation results for the problem described in Figure \ref{omega0}. We set the diffusive coefficient $k=1$ and compute the solution for two advection fields, constant on $\Omega$.
All cases were tested using a uniform mesh of 100 quadratic elements. We use the approximated RFB method described in sections 2 and 3 for computation. Convergence analysis for different  approximation orders and two P\'eclet numbers is illustrated on Figures \ref{Convergense1},\ref{Convergense2}. The convergence diagrams of the hp-FEM with equivalent approximation order are provided for reference. As can be seen, the approximated RFB method provides a much faster convergence speed compared to the hp-FEM when dealing with high P\'eclet numbers. The approximated RFB method has a clearly observable tendency to produce	 better approximation in $L_\infty$ norm when degrees of freedom is an odd number.
\begin{figure}
\center\includegraphics[width=0.65\textwidth]{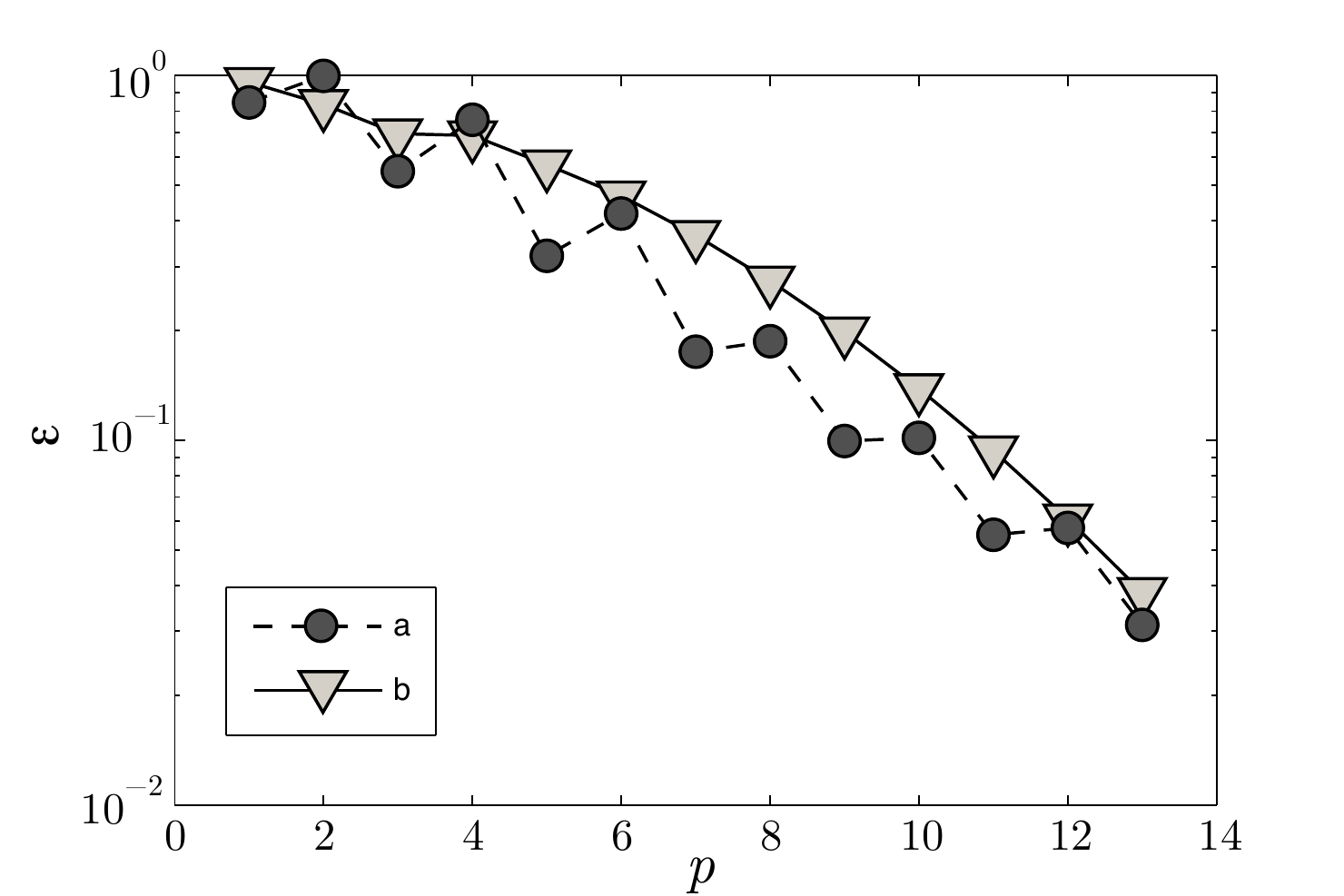}
\caption{$L_{\infty}$  error norm ($\varepsilon$) in nodes as a function of the approximation order $p$; a)approximated RFB method; b)hp-FEM. Mesh P\'eclet number $1.25\times10^2$.}
\label{Convergense1}
\end{figure}
\begin{figure}
\center\includegraphics[width=0.65\textwidth]{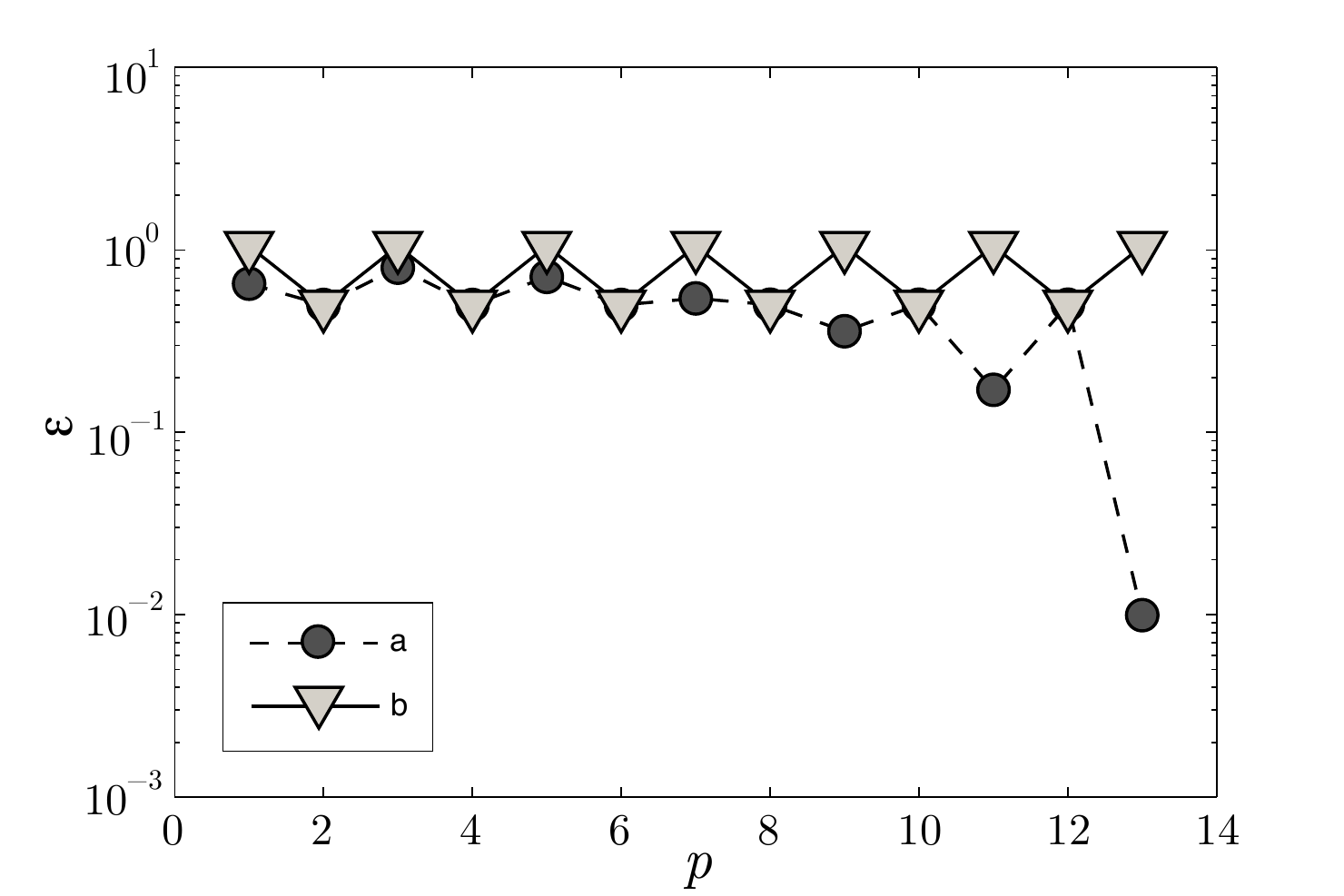}
\caption{$L_{\infty}$  error norm ($\varepsilon$) in nodes as a function of the approximation $p$; a)approximated RFB method; b)hp-FEM. Mesh P\'eclet number $1.25\times10^{14}$.}
\label{Convergense2}
\end{figure}

\subsection{A two-dimensional problem with non-symmetric boundary conditions}
\begin{figure}
\center
\includegraphics[width=0.5\textwidth]{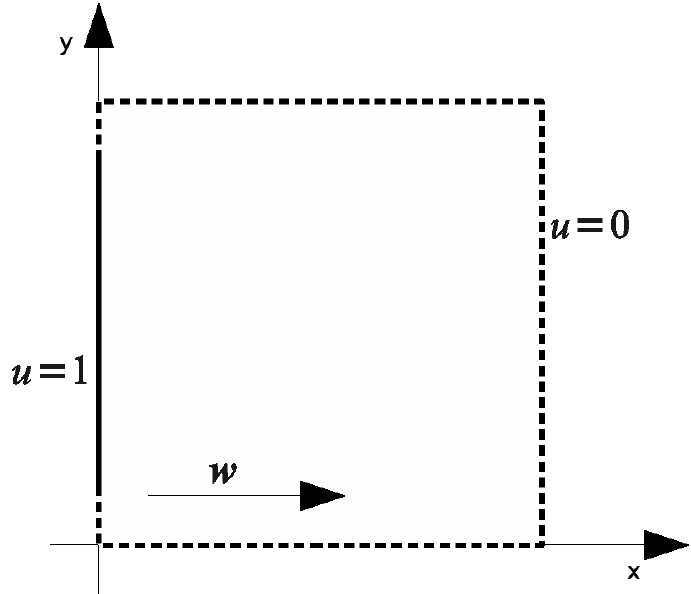}
\caption{2D domain with boundary condition}
\label{omega1}
\end{figure}
Let's consider problem described in Figure \ref{omega1}. As in the previous example we set $k=1$.
Since the stabilization effect is of prior interest, all cases were tested using 100 quadratic elements and uniform mesh. Computational results for the approximated RFB method described in sections 2,3 and for hp-FEM with $p=13$ are provided. We consider cases with different P\'eclet numbers (see Figures \ref{hres1},\ref{hres2},\ref{rres1},\ref{rres2}). 

It can be observed, the both methods exhibit a good stabilization effect when P\'eclet number is relatively small (Figures \ref{hres1},\ref{rres1}). However, this tendency is not kept when the P\'eclet number becomes higher. Indeed, hp-FEM shows non-physical oscillations when P\'eclet number reaches $10^5$ (Figure \ref{hres2}). In contrast, the approximated RFB method continues to generate stable results up to $Pe=10^{15}$, Figure \ref{rres2}.
The relationship between degrees of freedom and time, needed to compute the task, is illustrated in Figure \ref{time}. 

We have to note here, that due to the fact that the coefficients of \eqref{Lu} being constant in the whole domain $\Omega$, an obvious performance optimization for the developed methods is possible -- the basis enrichment $F_K(x)$ is computed once and used on every element. Since similar optimization could not be applied to hp-FEM  a large difference in computational performance is observe in Figure \ref{time}.
\begin{figure}
\centering
\includegraphics[width=0.6\textwidth]{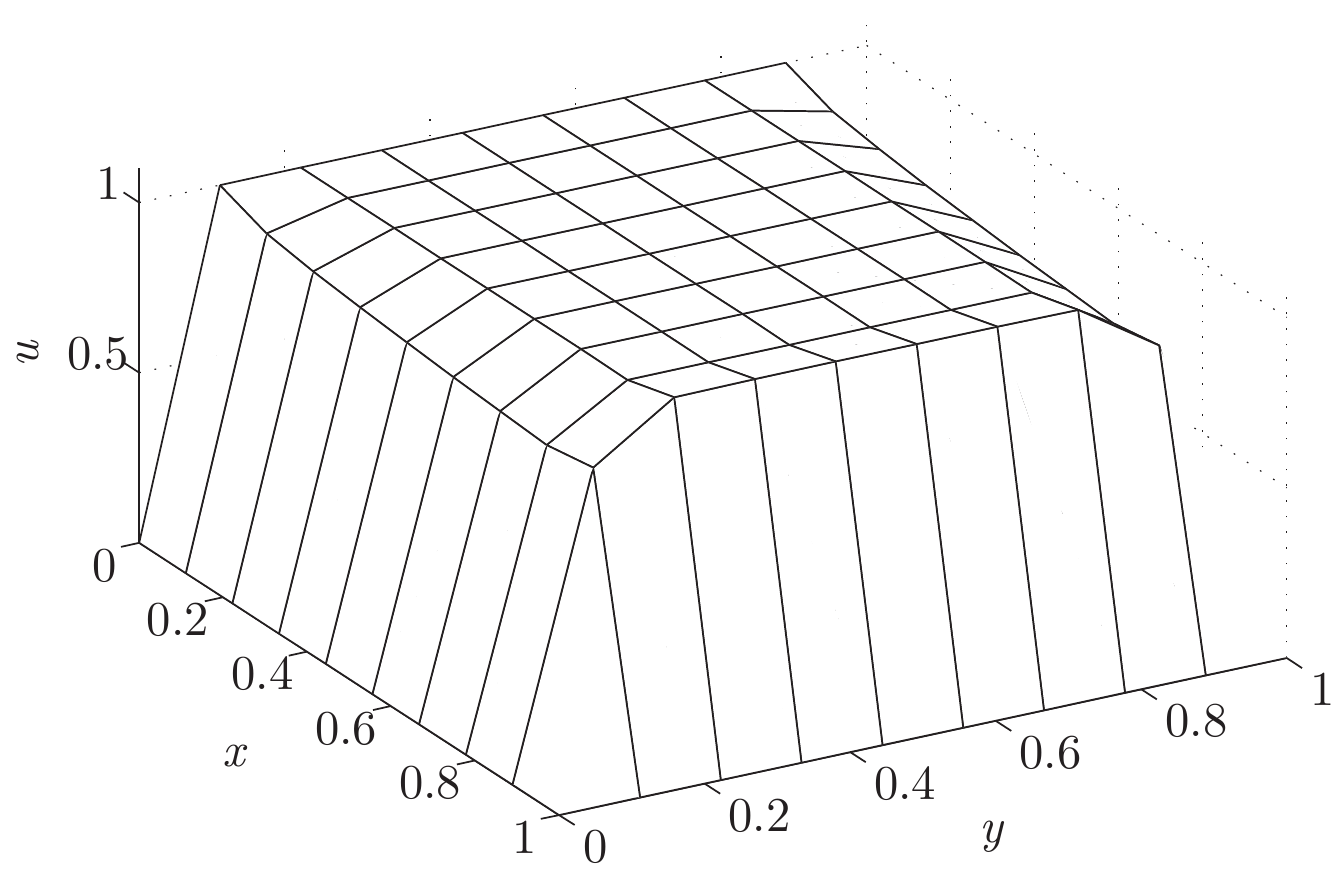}
\caption{Approximated solution, computed with hp-FEM, 100 elements, p=13, $Pe=10^3$.}
\label{hres1}
\end{figure}
\begin{figure}
\centering
\includegraphics[width=0.6\textwidth]{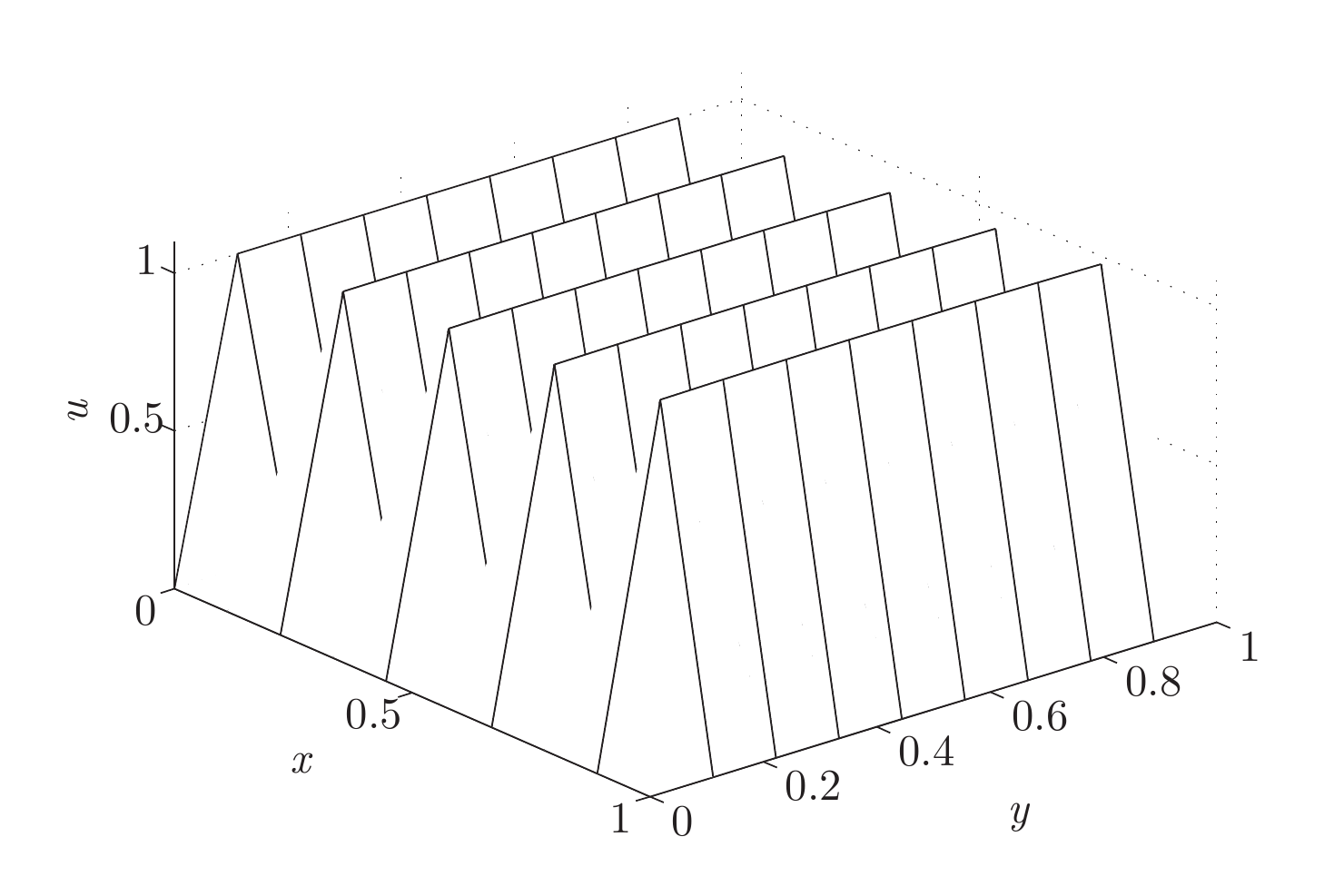}
\caption{Approximated solution, computed with hp-FEM, 100 elements, p=13, $Pe=10^5$.}
\label{hres2}
\end{figure}
\begin{figure}
\centering
\includegraphics[width=0.6\textwidth]{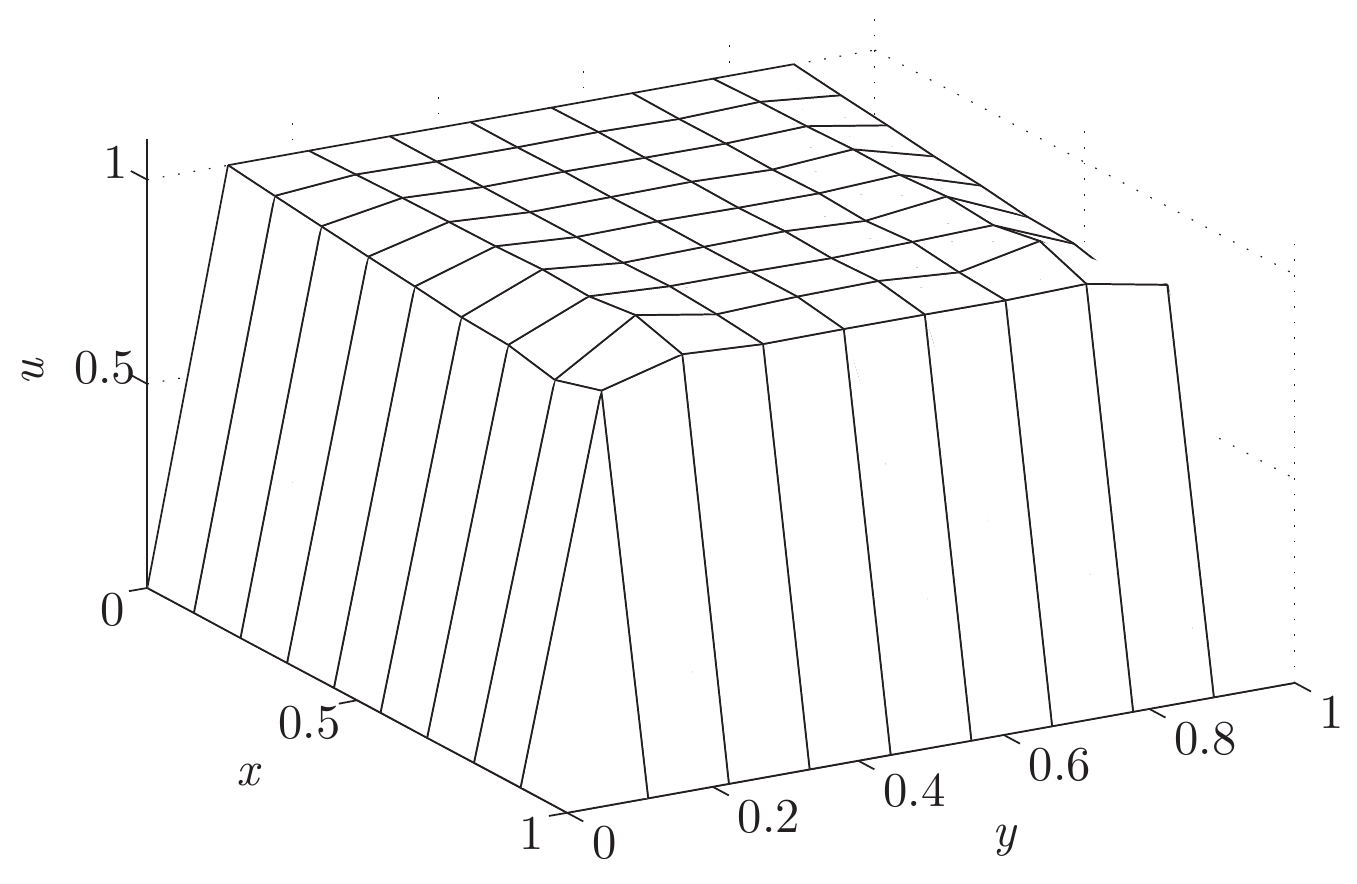}
\caption{Approximated solution, computed with approximated RFB method, 100 elements, p=13, $Pe=10^3$.}
\label{rres1}
\end{figure}
\begin{figure}
\centering
\includegraphics[width=0.6\textwidth]{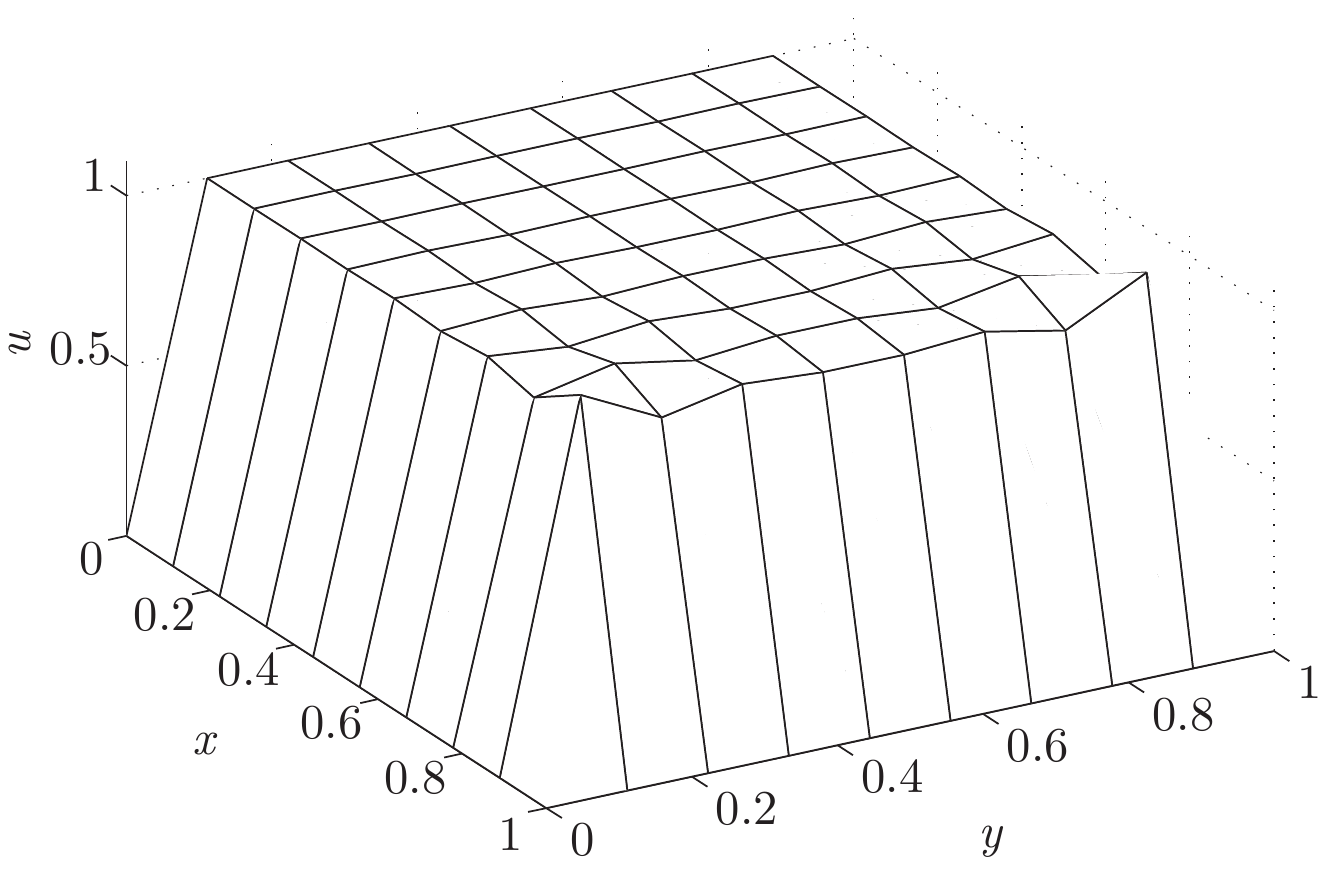}
\caption{Approximated solution computed with approximated RFB method, 100 elements, p=13, $Pe=10^{15}$.}
\label{rres2}
\end{figure}
\begin{figure}
\center
\includegraphics[width=0.6\textwidth]{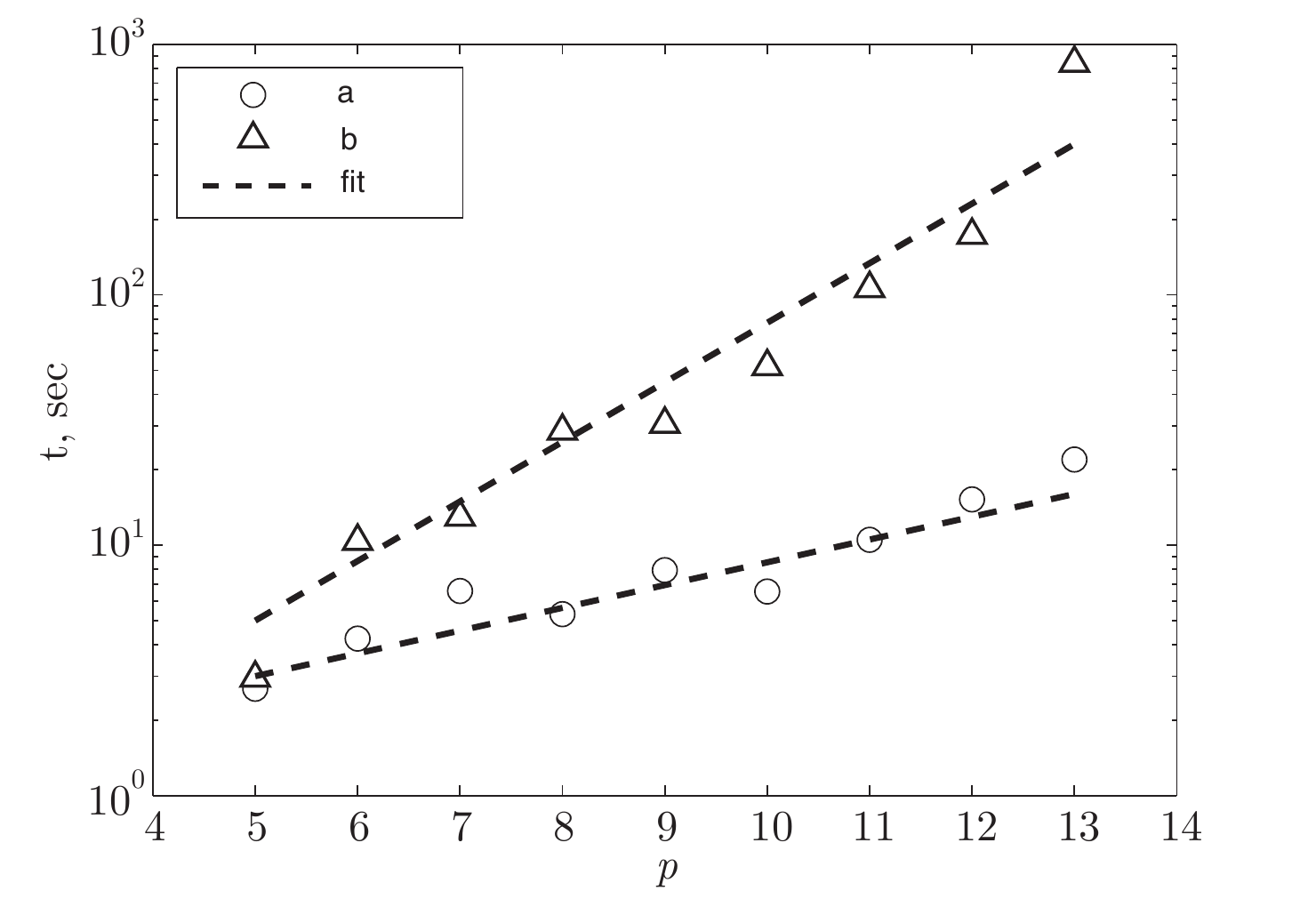}
\caption{The dependence of computational performance on approximation order for two methods: a)approximated RFB method; b)hp-FEM.}
\label{time}
\end{figure}
\section*{Conclusions}
A new strategy for implementing the approximated residual-free bubble method is proposed. It was found that the approach shows a strong stabilization effect when applied to the advection-diffusion problem with high P\'eclet numbers. Using the proposed method, oscillation-free numerical solutions of the target problem were achieved for mesh P\'eclet numbers up to $10^{14}$.  The proposed discrete scheme shows high computational performance, and is parallelization friendly.

Furthermore, the strategy finds the approximated solution of a PDE in the piece-polynomial space using a basis similar to that of hp-FEM. From this point of view, one may refer to the discussed strategy as a substitution to hp-FEM. When comparing the numerical results to those of hp-FEM, a better convergence and computational performance is observed.
The use of odd number for degrees of freedom is recommended.
\section*{References}
\bibliographystyle{plain}

\bibliography{bibliography} 
\end{document}